\documentclass[reqno,twoside]{article}
\usepackage{amssymb,amsmath,amsthm}
\usepackage[paperwidth=200mm,paperheight=260mm,left=28.4mm,right=28.4mm,top=25.6mm,bottom=21.9mm,headsep=5mm]{geometry}
\usepackage{mathrsfs}
\usepackage{cite}

\newtheorem{thm}{Theorem}[section]

\newtheorem{prop}[thm]{Proposition}

\newtheorem{rem}[thm]{Remark}
\numberwithin{equation}{section} \allowdisplaybreaks
\newcommand{\norm}[1]{\left\|#1\right\|}
\newcommand{\nnorm}[1]{\left|\!\left|\!\left|#1\right|\!\right|\!\right|}
\newcommand{\abs}[1]{\left|#1\right|}

\newcommand{\Real}{\mathbb R}

\newcommand{\no}{\nonumber}
\newcommand{\eps}{\varepsilon}

\newcommand{\T}{\mathscr{T}}

\newcommand{\F}{\mathscr{F}}

\renewcommand{\Re}{\mathrm{Re}}
\newcommand{\ls}{\leqslant}
\newcommand{\gs}{\geqslant}
\newcommand{\dif}{\mathrm{d} }

\newenvironment{pf}[1][Proof]{{\par \emph{#1.}}\;}{\hfill $\Box$\par}
\begin{document}
\pagestyle{myheadings} \markboth{\hfill \textrm{C. C. HAO, L. HSIAO
AND H. L. Li} \hfill}{\hfill \textrm{WELL-POSEDNESS FOR ROTATING
GPE} \hfill}
\title{\Large\textbf{Global well-posedness for the Gross-Pitaevskii
equation with an angular momentum rotational term}}%
\author{\large Chengchun Hao${}^{1,}$\thanks{Correspondence to: Chengchun Hao, Institute of Mathematics,
Academy of Mathematics \& Systems Science, CAS, Beijing 100190,
People's Republic of China.} ${}^{,}$\thanks{E-mail:
hcc@amss.ac.cn\newline\newline Contract/grant sponsor: National
Natural Science Foundation of China\newline Contract/grant sponsor:
Chinese Academy of Sciences\newline Contract/grant sponsor: Beijing
Nova Program}
, Ling Hsiao${}^1$ and Hai-Liang Li${}^2$\\[1mm]
\footnotesize\textit{${}^1$Institute of Mathematics, Academy of
Mathematics \& Systems
Science, CAS,}\\
\footnotesize \textit{Beijing 100190, People's Republic of China}\\
\footnotesize \textit{${}^2$Department of Mathematics, Capital
Normal University},\\
\footnotesize \textit{Beijing 100037, People's Republic of China}}
\date{}
 \maketitle
 \thispagestyle{empty}

\begin{center}\textbf{SUMMARY}\end{center}

\noindent In this paper, we establish the global well-posedness of
the Cauchy problem for the Gross-Pitaevskii equation with an
rotational angular momentum term in the space $\Real^2$.\\



\section{Introduction}

The Gross-Pitaevskii equation (GPE), derived independently by Gross
\cite{Gro61} and Pitaevskii \cite{Pit61}, arises in various models
of nonlinear physical phenomena. This is a Schr\"odinger-type
equation with an external field potential $V_{ext}(t,x)$ and a local
cubic nonlinearity:
\begin{align}\label{cgpe}
    i\hbar\partial_t u+\frac{\hbar^2}{2m}\Delta u=V_{ext}u+\beta
    \abs{u}^2u.
\end{align}

The GPE \eqref{cgpe} in physical dimensions ($2$ and $3$ dimensions)
is used in the meanfield quantum theory of Bose-Einstein condensate
(BEC) formed by ultracold bosonic coherent atomic ensembles.
Recently, several research groups \cite{HMEWC,MCWD00,MCBD01,MAHHWC}
have produced quantized vortices in trapped BECs, and a typical
method they used is to impose a laser beam on the magnetic trap to
create a harmonic anisotropic rotating trapping potential. The
properties of BEC in a rotational frame at temperature $T$ being
much smaller than the critical condensation temperature $T_c$
\cite{LL77} are well described by the macroscopic wave function
$u(t,x)$, whose evolution is governed by a self-consistent, mean
field nonlinear Schr\"odinger equation (NLS) in a rotational frame,
also known as the Gross-Pitaevskii equation with an angular momentum
rotation term:
\begin{align}\label{rgpe}
    i\hbar\partial_t u+\frac{\hbar^2}{2m}\Delta u=V(x)u+NU_0
    \abs{u}^2u-\Omega L_z u, \; x\in\Real^3,\, t\gs 0,
\end{align}
where the wave function $u(t,x)$ corresponds to a condensate state,
$m$ is the atomic mass, $\hbar$ is the Planck constant, $N$ is the
number of atoms in the condensate, $\Omega$ is the angular velocity
of the rotating laser beam, and $V(x)$ is an external trapping
potential. When a harmonic trap potential is considered,
$V(x)=\frac{m}{2}\left(\omega_{1}^2x_1^2+\omega_{2}^2x_2^2+\omega_{3}^2x_3^2\right)$
with $\omega_{1}$, $\omega_{2}$ and $\omega_{3}$ being the trap
frequencies in the $x_1$-, $x_2$- and $x_3$-direction, respectively.
The local nonlinearity term $NU_0\abs{u}^2u$ arises from an
assumption about the delta-shape interatomic potential. $U_0=4\pi
\hbar^2 a_s/m$ describes the interaction between atoms in the
condensate with $a_s$ (positive for repulsive interaction and
negative for attractive interaction) the $s$-wave scattering length,
and $L_z=-i\hbar(x_1\partial_{x_2}-x_2\partial_{x_1})$ is the third
component of the angular momentum $L=x\times P$ with the momentum
operator $P=-i\hbar\nabla$.

After normalization, proper nondimensionalization and dimension
reduction in certain limiting trapping frequency regime \cite{PS03},
it turns to be the dimensionless GPE in $d$-dimensions ($d=2,3$):
\begin{align}\label{dlgpe}
    iu_t+\frac{1}{2}\Delta u=V_d(x)u+\beta_d\abs{u}^2u-\Omega L_z
    u,\;x\in\Real^d,\, t>0,
\end{align}
where $L_z=i(x_1\partial_{x_2}-x_2\partial_{x_1})$ and
\begin{align}\label{V}
    \beta_d=\left\{\begin{array}{l}\beta\sqrt{\gamma_3/2\pi},\\
    \beta,\end{array}\right.\quad V_d(x)=\left\{\begin{array}{ll}
    (\gamma_1^2x_1^2+\gamma_2^2x_2^2)/2, \;&d=2,\\
    (\gamma_1^2x_1^2+\gamma_2^2x_2^2+\gamma_3^2x_3^2)/2,
    &d=3,\end{array}\right.
\end{align}
with $\gamma_1>0$, $\gamma_2>0$ and $\gamma_3>0$ constants,
$\beta=\frac{4\pi a_s N}{a_0}$, $a_0=\sqrt{\frac{\hbar}{m\omega_m}}$
and $\omega_m=\min\{\omega_1,\omega_2,\omega_3\}$.

In general, it is a rather complicated process about the dynamics of
solutions (in particular, vortex) for GPE \eqref{rgpe} under the
interaction of trapping frequencies and angular rotating motion. The
recent numerical simulation of GPE \eqref{rgpe} for different choice
of trap frequencies $(\gamma_1,\gamma_2)$ can help us to understand
the complicated dynamical phenomena caused by the angular rotating
and spatial high frequency motion. The case of different frequency
$\gamma_1\neq\gamma_2$ gives much complicated behavior and thus is
rather difficult to be studied rigorously \cite{BDZ06,BWM05}. To our
knowledge, the equation \eqref{rgpe} has been only investigated for
some specific cases by numerical simulation. Therefore, to develop
methods for constructing analytical solutions of the GPE
\eqref{cgpe} or some specific cases is the first step in order to
understand the dynamics caused by the trapping and rotation.

To begin with, we first consider the case $\gamma_1=\gamma_2=\omega$ which
means the spatial isotropic motion and focus on the Cauchy problem
of the Gross-Pitaevskii equation with an angular momentum rotational
term in two dimensions
\begin{align}\label{gpe}
    &iu_t+\frac{1}{2}\Delta u
    =\frac{\omega^2}{2}\abs{x}^2 u
    +\beta\abs{u}^{2\sigma}u-\omega L_z u,\;
    x\in\Real^2,\, t\gs 0, \\
    &u(0,x)=u_0(x),\, x\in\Real^2,\label{data}
\end{align}
where the wave function $u=u(t,x):
[0,\infty)\times\Real^2\to\mathbb{C}$ corresponds to a condensate
state, $\Delta$ is the Laplace operator on $\Real^2$, $\omega>0$,
$\beta>0$ and $\sigma\in [1/\omega,\infty)$ are constants, and
$L_z=-i(x_1\partial_{x_2}-x_2\partial_{x_1})=i(x_2\partial_{x_1}-x_1\partial_{x_2})$
is the dimensionless angular momentum rotational term
\cite{BDZ06,BWM05}. We assume that the initial value
\begin{align}\label{data.space}
    u_0(x)\in \Sigma:=\{u\in H^1(\Real^2):\, \abs{x}u\in
    L^2(\Real^2)\},
\end{align}
with the norm
\begin{align*}
    \norm{u}_\Sigma =\norm{u}_{H^1}+\norm{\abs{x}u}_{L^2}.
\end{align*}

It is clear that \eqref{gpe} is a special case of \eqref{dlgpe} with
a spatial isotropic trapping frequency (i.e. $\gamma_1=\gamma_2$)
and $\Omega=\omega$ in $2$-dimensions. Note that for
multi-dimensional GPE~\eqref{rgpe}, nothing is known about the exact
integration except for the case when $\gamma_1=\gamma_2$
($=\gamma_3$ for 3D) without the angular momentum rotational term,
namely, $\Omega=0$ considered in \cite{Car02,Car03}.

There are three ingredients that play important roles in the proof
of our result. The first involves the solution of the Cauchy problem
to the linear equation
\begin{align}\label{lin.eqn}
\begin{split}
    &iu_t+\frac{1}{2}\Delta u=\frac{\omega^2}{2}\abs{x}^2 u-\omega L_z
    u,\;
    x\in\Real^2,\, t\gs 0, \\
    &u(0,x)=u_0(x),\, x\in\Real^2,
\end{split}
\end{align}
which is significant for investigating the properties of the
evolution operator corresponding to the linear operator
$i\partial_t+\frac{1}{2}\Delta-\frac{\omega^2}{2}\abs{x}^2+\omega
L_z$. The second one is to obtain the Strichartz estimates for the
foregoing linear operator. The last one is that there exist two
Galilean operators $J(t)$ and $H(t)$ (as blow) which can commute
with the linear operator and can be viewed as the substitute of
$\nabla$ and $x$ respectively in the non-potential case.

Now we state the main result of this paper.

\begin{thm}\label{thm}
Let $u_0\in\Sigma$ and $\rho\in[2,\infty)$. Then, there exists a
unique solution $u(t,x)$ to the Cauchy problem
\eqref{gpe}--\eqref{data}. And the solution satisfies, for any
$T\in(0,\infty)$, that
\begin{align*}
    u(t,x), J(t)u(t,x), H(t)u(t,x)\in \mathcal{C}(\Real;
    L^2(\Real^2))\cap L^{\gamma(\rho)}(0,T; L^\rho(\Real^2)),\;
    \forall t\in [0,T],
\end{align*}
where $\frac{1}{\gamma(\rho)}=\frac{1}{2}-\frac{1}{\rho}$, $J(t)$
and $H(t)$ are defined as below as in \eqref{J} and \eqref{H},
respectively.
\end{thm}

\begin{rem}
Since the GPE \eqref{gpe} (or \eqref{dlgpe}) in a rotational frame
is time reversible and time transverse invariant, the above result
is also valid for the case when $t<0$.
\end{rem}

The paper is organized as follows.  In Sec.~\ref{Sec.Str}, the
evolution operator of the linear equation and the Strichartz
estimates about the former operator are first established.
Sec.~\ref{Sec.law} is devoted to the derivation of some conservation
identities such as the mass, the energy, the angular momentum
expectation, the pseudo-conformal conservation laws in the whole
space $\Real^2$ for \eqref{gpe}--\eqref{data}. Finally, the
nonlinear estimates and the proof of Theorem~\ref{thm} are obtained
in Sec.~\ref{Sec.pf}.

\section{The Strichartz estimates and some main
operators}\label{Sec.Str}

Let $u(t)$ be the solution of the linear equation \eqref{lin.eqn},
then by a computation, it can be expressed as
\begin{align}\label{def.S}
    u(t)=S(t)u_0=\frac{\omega}{2\pi i\sin (\omega
    t)}\int_{\Real^2}e^{i\omega\left(\frac{\abs{x-y}^2}{2}\cot(\omega t)-x^\perp\cdot
    y\right)} u_0(y)dy,
\end{align}
where $x^\perp :=(-x_2,x_1)$ and $S(t)$ is the evolution operator
which can be formally written as
$S(t):=e^{i\frac{t}{2}(\nabla-i\omega x^\perp)^2}$.  This formula,
which can be deduced from the three-dimensional nonlinear
Schr\"odinger equation with a magnetic field discussed in
\cite[Sec.9.1]{Caz03} or \cite{AHS78}, defines a operator $S(t)$,
unitary on $L^2$. Note that this formula is valid only for small time, due to the singularity formation for the fundamental solution. For this nonlinear Schr\"odinger equation,
Strichartz estimates are available. These estimates, mixed
time-space estimates, are exactly the same as for
$S_0(t)=e^{\frac{i}{2}t\Delta}$. Recall the main properties from
which such estimates stem. The operator $S(t)$ is unitary on $L^2$,
$\norm{S(t)}_{L^2\to L^2}=1$. In fact, by the Plancherel theorem, we
have for any $\phi\in L^2$
\begin{align*}
    \norm{S(t)\phi}_2=&\frac{\omega}{2\pi\abs{\sin \omega
    t}}\norm{\int_{\Real^2}e^{i\omega\left(\frac{\abs{x-y}^2}{2}\cot(\omega t)-x^\perp\cdot
    y\right)} \phi(y)dy}_2\\
    =&\frac{\omega }{2\pi\abs{\sin \omega
    t}}\norm{e^{\frac{i\omega\abs{x}^2\cot (\omega t)}{2}}\int_{\Real^2}e^{i\omega\left(x\cot (\omega t)-x^\perp\right)\cdot
    y}e^{i\omega\frac{\abs{y}^2}{2}\cot(\omega t)} \phi(y)dy}_2\\
    =&\frac{\omega }{\abs{\sin \omega
    t}}\norm{\left[\F^{-1}\left(e^{i\omega\frac{\abs{y}^2}{2}\cot(\omega t)}
    \phi(y)\right)\right]\left(\omega\left(x\cot (\omega t)-x^\perp\right)\right)
    }_2\\
    =&\norm{\F^{-1}\left(e^{i\omega\frac{\abs{y}^2}{2}\cot(\omega t)}
    \phi(y)\right)}_2
    =\norm{e^{i\omega\frac{\abs{y}^2}{2}\cot(\omega t)}
    \phi(y)}_2\\
    =&\norm{\phi}_2.
\end{align*}
And for $0<t\ls \frac{\pi}{2\omega}$, the operator is dispersive,
with $\norm{S(t)}_{L^1\to L^\infty}\ls \frac{1}{4}\abs{t}^{-1}$,
since $\abs{\sin t}\gs \frac{2}{\pi}\abs{t}$ for
$\abs{t}\ls\frac{\pi}{2}$. Thus, we can obtain similar Strichartz
estimates to the linear Schr\"odinger operator
$e^{i\frac{t}{2}\Delta}$ by the standard methods (c.f. \cite{KT98})
provided that only finite time intervals are involved
(c.f.\cite{Car03}).

\begin{prop}\label{pr.Str}
Let $I$ be an interval contained in $[0,\pi/2\omega]$. Then, it
holds that

\mbox{\rm (1)} For any admissible pair $(\gamma(p),p)$ (that is,
$1/\gamma(p)=1/2-1/p$ for $2\ls p<\infty$), there exists $C_p$ such
that for any $\phi\in L^2$
\begin{align}\label{Str.1}
    \norm{S(t)\phi}_{L^{\gamma(p)}(I;L^p)}\ls C_p\norm{\phi}_{L^2}.
\end{align}

\mbox{\rm (2)} For any admissible pairs $(\gamma(p_1),p_1)$ and
$(\gamma(p_2),p_2)$, there exists $C_{p_1,p_2}$ such that
\begin{align}\label{Str.2}
    \norm{\int_{I\cap\{s\ls t\}}S(t-s)F(s)\dif
    s}_{L^{\gamma(p_1)}(I;L^{p_1})}\ls
    C_{p_1,p_2}\norm{F}_{L^{\gamma(p_2)'}(I;L^{p'_2})}.
\end{align}
The above constants are independent of $I\subset[0,\pi/2\omega]$.
\end{prop}

The integral equation reads
\begin{align}\label{e.int}
    u(t)=S(t)u_0-i\beta \int_0^t S(t-s)\abs{u}^{2\sigma}
    u(s) \dif s.
\end{align}

Since the initial data belong to the space $\Sigma$, we naturally
need the estimates of $\nabla S(t) \phi$ and $xS(t)\phi$. In fact,
from \eqref{def.S}, we can compute and obtain that
\begin{align*}
    \nabla S(t)\phi=i\omega x\cot (\omega t) S(t)\phi-i\omega \cot(\omega t)
    S(t)(x\phi)+i\omega S(t)(x^\perp \phi),
\end{align*}
and
\begin{align*}
    &S(t)\nabla^\perp \phi=i\omega x^\perp\cot(\omega t)
    S(t)\phi-i\omega \cot (\omega t) S(t)(x^\perp \phi)-i\omega
    xS(t)\phi,\\
    &S(t)\nabla \phi=i\omega x\cot(\omega t)
    S(t)\phi-i\omega \cot (\omega t) S(t)(x \phi)+i\omega
    x^\perp S(t)\phi,
\end{align*}
which yield to
\begin{align}
    \nabla S(t)\phi=&\cos(\omega t) S(t)(\cos(\omega t)\nabla- \sin(\omega
    t)\nabla^\perp)\phi\no\\
    &-i\omega \sin (\omega t) S(t)\left[(\cos(\omega t)x-\sin(\omega t)x^\perp
    )\phi\right],\label{e1}\\
    xS(t)\phi=&\cos(\omega t) S(t)\left[(\cos(\omega t)x-\sin(\omega t)x^\perp
    )\phi\right]\no\\
    &-\frac{i}{\omega}\sin (\omega t) S(t)(\cos(\omega t)\nabla- \sin(\omega
    t)\nabla^\perp)\phi.\label{e2}
\end{align}

Thus, we have
\begin{align*}
    &S(t)(-i\nabla)\phi\no\\
    =&\left[\omega \sin (\omega t)(\cos(\omega t)x+\sin(\omega t)x^\perp)
    -i\cos(\omega t)(\cos(\omega t)\nabla+ \sin(\omega
    t)\nabla^\perp)\right]S(t)\phi,
\end{align*}
and
\begin{align*}
    &S(t)\omega x\phi\no\\
    =&\left[\omega \cos (\omega t)(\cos(\omega t)x+\sin(\omega t)x^\perp)
    +i\sin(\omega t)(\cos(\omega t)\nabla+ \sin(\omega
    t)\nabla^\perp)\right]S(t)\phi.
\end{align*}

 For convenience of computations, we denote
\begin{align}\label{J}
    J(t)=\omega \sin (\omega t)(\cos(\omega t)x+\sin(\omega t)x^\perp)
    -i\cos(\omega t)(\cos(\omega t)\nabla+ \sin(\omega
    t)\nabla^\perp),
\end{align}
and the corresponding ``orthogonal'' operator
\begin{align}\label{H}
    H(t)=\omega \cos (\omega t)(\cos(\omega t)x+\sin(\omega t)x^\perp)
    +i\sin(\omega t)(\cos(\omega t)\nabla+ \sin(\omega
    t)\nabla^\perp),
\end{align}
which will appear in the pseudo-conformal conservation law and play
a crucial role in the nonlinear estimates.

Thus, we get
\begin{align*}
    J(t)=S(t)(-i\nabla)S(-t), \quad H(t)=S(t)\omega xS(-t).
\end{align*}

By computation, we can obtain the following commutation relation
\begin{align}\label{commu.rel}
\begin{split}
    &\left[J(t),i\partial_t+\frac{1}{2}\Delta -\frac{\omega^2}{2}\abs{x}^2 +\omega
    L_z\right]=0,\\
    &\left[H(t),i\partial_t+\frac{1}{2}\Delta -\frac{\omega^2}{2}\abs{x}^2 +\omega
    L_z\right]
    =0.
\end{split}
\end{align}

In addition, denote $M(t)=e^{-i\omega \frac{\abs{x}^2}{2}\tan(\omega
t)}$ and $Q(t)=e^{i\omega \frac{\abs{x}^2}{2}\cot(\omega t)}$, then
\begin{align}\label{JH}
\begin{split}
    J(t)=&-i\cos(\omega t) M(t)(\cos(\omega t) \nabla+\sin (\omega t)
    \nabla^\perp)M(-t),\\
    H(t)=&i\sin(\omega t) Q(t)(\cos(\omega t) \nabla+\sin (\omega t)
    \nabla^\perp)Q(-t).
\end{split}
\end{align}

\section{The conserved quantities}\label{Sec.law}

\begin{prop}\label{pr.cons}
Let $u$ be a solution of the equation \eqref{gpe} with the initial
data $\phi\in\Sigma(\Real^2)$. Then, we have the following conserved
quantities for all $t\gs 0$:

 \mbox{\rm (1)}
The $L^2$-norm:
\begin{align}\label{e.mass}
    \norm{u(t)}_{2}=\norm{u_0}_{2}.
\end{align}

 \mbox{\rm (2)}
The energy for the non-rotating part:
\begin{align}\label{e.non}
    E_0(u)=\frac{1}{2}\norm{\nabla
    u}_2^2+\frac{\omega^2}{2}\norm{x
    u}_2^2+\frac{\beta}{\sigma+1}\norm{u}_{2\sigma+2}^{2\sigma+2}
    =E_0(u_0).
\end{align}

\mbox{\rm (3)} The angular momentum expectation:
\begin{align}\label{e.mom}
    \langle L_z\rangle(t)=\int_{\Real^2}\bar{u}L_z u \dif x =\langle
    L_z\rangle(0).
\end{align}
\end{prop}

\begin{pf}
For convenience, we introduce
\begin{align*}
    eq(u):=iu_t+\frac{1}{2}\Delta u-\frac{\omega^2}{2}\abs{x}^2 u-\beta\abs{u}^{2\sigma}u+\omega L_z
    u.
\end{align*}

It is clear that \eqref{e.mass} holds by applying the $L^2$-inner
product between $eq(u)$ and $\bar{u}$, and then taking the imaginary
part of the resulting equation.

Since we can use the identity \eqref{e.mom} in the proof of
\eqref{e.non}, we derive \eqref{e.mom} first.  Differentiating
$\langle L_z\rangle(t)$ with respect to $t$, and integrating by
parts, we have
\begin{align*}
    \frac{d\langle L_z\rangle(t)}{dt}
    =&i\int_{\Real^2}\left[\bar{u}_t(x_2\partial_{x_1}u-x_1\partial_{x_2}u)
    +\bar{u}(x_2\partial_{x_1}u_t-x_1\partial_{x_2}u_t)\right] \dif
    x\\
    =&\int_{\Real^2}\left[\overline{-iu_t}(x_2\partial_{x_1}u-x_1\partial_{x_2}u)
    -iu_t(x_2\partial_{x_1}\bar{u}-x_1\partial_{x_2}\bar{u})\right] \dif
    x\\
    =&\int_{\Real^2}\left[\left(\frac{1}{2}\Delta \bar{u}-\frac{\omega^2}{2}\abs{x}^2 \bar{u}
    -\beta\abs{u}^{2\sigma}\bar{u}+\omega L_z
    \bar{u}\right)(x_2\partial_{x_1}u-x_1\partial_{x_2}u)\right.\\
    &\qquad\left.+\left(\frac{1}{2}\Delta u-\frac{\omega^2}{2}\abs{x}^2 u-\beta\abs{u}^{2\sigma}u+\omega L_z
    u\right)(x_2\partial_{x_1}\bar{u}-x_1\partial_{x_2}\bar{u})\right] \dif
    x\\
    =&\int_{\Real^2}\Re(\Delta
    u(x_2\partial_{x_1}\bar{u}-x_1\partial_{x_2}\bar{u}))
    -\frac{\omega^2}{2}(\abs{x}^2(x_2\partial_{x_1}\abs{u}^2-x_1\partial_{x_2}\abs{u}^2))\\
    &\qquad
    -\beta\Re(\abs{u}^{2\sigma}(x_2\partial_{x_1}\abs{u}^2-x_1\partial_{x_2}\abs{u}^2))\dif
    x\\
    =&\frac{\omega^2}{2}\int_{\Real^2}(2x_1x_2\abs{u}^2-2x_2x_1\abs{u}^2))\dif
    x\\
    =&0,
\end{align*}
which yields the desired identity \eqref{e.mom}.

Next, we prove the energy conservation law for the non-rotating part
\eqref{e.non}. We consider
\begin{align*}
    \Re (eq(u),u_t)=0,
\end{align*}
where $(\cdot,\cdot)$ denotes the $L^2$-inner product. From the
above, we can get
\begin{align*}
    \int_{\Real^2}\left[\frac{1}{2}\partial_t\abs{\nabla
    u}^2+\frac{\omega^2}{2}\partial_t\abs{xu}^2+\frac{\beta}{\sigma+1}\partial_t\abs{u}^{2\sigma+2}
    +\frac{\omega}{2}\partial_t(\bar{u}L_zu)\right]\dif x=0,
\end{align*}
which implies the identity \eqref{e.non} with the help of
\eqref{e.mom}.
\end{pf}

\begin{rem} For the equation \eqref{gpe}, the pseudo-conformal type conservation laws are also valid:
\begin{align*}
\begin{split}
   &\norm{H(t) u}_2^2 +\frac{2\beta \sin^2(\omega t)}{\sigma+1}\norm{u}_{2\sigma+2}^{2\sigma+2}
   +\frac{2\beta[\sigma\omega-1]}{\sigma+1}\int_0^t\sin 2\omega s
    \norm{u(s)}_{2\sigma+2}^{2\sigma+2}\dif s\\
    =&\omega^2\norm{xu_0}_2^2,
\end{split}
\end{align*}
and
\begin{align*}
\begin{split}
   \norm{J(t) u}_2^2 +\frac{2\beta \cos^2(\omega t)}{\sigma+1}\norm{u}_{2\sigma+2}^{2\sigma+2}
    =&\norm{\nabla u_0}_2^2+\frac{2\beta}{\sigma+1}\norm{u_0}_{2\sigma+2}^{2\sigma+2}\\
    &+\frac{2\beta[\sigma\omega-1]}{\sigma+1}\int_0^t\sin 2\omega s
    \norm{u(s)}_{2\sigma+2}^{2\sigma+2}\dif s.
\end{split}
\end{align*}
\end{rem}

\section{Nonlinear estimates and the proof of
Theorem~\ref{thm}}\label{Sec.pf}

By computation, we can get, with the help of \eqref{JH}, that
\begin{align*}
    J(t)\abs{u}^{2\sigma}u=(\sigma
    +1)\abs{u}^{2\sigma}J(t)u-\sigma\abs{u}^{2\sigma-2}u^2\overline{J(t)u},
\end{align*}
which implies, in view of
$\frac{1}{\rho'}+\eps=\frac{2\sigma}{q}+\frac{1}{\rho}$ with
$0<\eps<\frac{1}{\rho}$, that
\begin{align*}
    \norm{J(t)\abs{u}^{2\sigma}u}_{L^{\left(\frac{\rho}{1-\rho\eps}\right)'}}\ls
    C\norm{u}_{L^q}^{2\sigma}\norm{J(t)u}_{L^\rho}.
\end{align*}
From the Sobolev embedding theorem and the H\"older inequality, it
yields
\begin{align*}
    &\norm{J(t)\abs{u}^{2\sigma}u}_{L^{\gamma\left(\frac{\rho}{1-\rho\eps}\right)'}(0,T;
    L^{\left(\frac{\rho}{1-\rho\eps}\right)'})}\\
    \ls &
    CT^{1-\eps-\frac{2}{\gamma(\rho)}}\norm{u}_{L^\infty(0,T; H^1)}^{2\sigma}
    \norm{J(t)u}_{L^{\gamma(\rho)}(0,T; L^\rho)}.
\end{align*}

Similarly, we have
\begin{align*}
    &\norm{H(t)\abs{u}^{2\sigma}u}_{L^{\gamma\left(\frac{\rho}{1-\rho\eps}\right)'}(0,T;
    L^{\left(\frac{\rho}{1-\rho\eps}\right)'})}\\
    \ls &
    CT^{1-\eps-\frac{2}{\gamma(\rho)}}\norm{u}_{L^\infty(0,T; H^1)}^{2\sigma}
    \norm{H(t)u}_{L^{\gamma(\rho)}(0,T; L^\rho)},
\end{align*}
and
\begin{align*}
    \norm{\abs{u}^{2\sigma}u}_{L^{\gamma\left(\frac{\rho}{1-\rho\eps}\right)'}(0,T;
    L^{\left(\frac{\rho}{1-\rho\eps}\right)'})}\ls
    CT^{1-\eps-\frac{2}{\gamma(\rho)}}\norm{u}_{L^\infty(0,T; H^1)}^{2\sigma}
    \norm{u}_{L^{\gamma(\rho)}(0,T; L^\rho)}.
\end{align*}

For convenience, we denote
\begin{align*}
    \nnorm{u}_A:=\norm{u}_A+\norm{J(t)u}_A+\norm{H(t)u}_A,
\end{align*}
where $A$ denotes a normalized space. Thus, we have
\begin{align}\label{non.est}
    &\nnorm{\abs{u}^{2\sigma}u}_{L^{\gamma\left(\frac{\rho}{1-\rho\eps}\right)'}(0,T;
    L^{\left(\frac{\rho}{1-\rho\eps}\right)'})}\no\\
    \ls &
    CT^{1-\eps-\frac{2}{\gamma(\rho)}}\norm{u}_{L^\infty(0,T; H^1)}^{2\sigma}
    \nnorm{u}_{L^{\gamma(\rho)}(0,T; L^\rho)}.
\end{align}

For any  $\rho\in [2,\infty)$ and $M\gs 2C\norm{u_0}_\Sigma$, define
the workspace $(\mathcal{D},d)$ as
\begin{align*}
    \mathcal{D}:=\{u:\, \nnorm{u}_{L^\infty(0,T;L^2)\cap L^{\gamma(\rho)}(0,T; L^\rho)}\ls M\},
\end{align*}
with the distance
\begin{align*}
    d(u,v)=\nnorm{u-v}_{L^{\gamma(\rho)}(0,T; L^\rho)}.
\end{align*}
It is clear that $(\mathcal{D},d)$ is a Banach space. Let us
consider the mapping $\T :(\mathcal{D},d)\to (\mathcal{D},d)$
defined by
\begin{align*}
    \T: u(t)\mapsto S(t)u_0-i\beta \int_0^t S(t-s)\abs{u}^{2\sigma}
    u(s) \dif s.
\end{align*}

For $u\in (\mathcal{D},d)$, by the commutation relation
\eqref{commu.rel}, Proposition~\ref{pr.Str} and the nonlinear
estimate \eqref{non.est}, we obtain
\begin{align}\label{non1}
    \nnorm{\T u}_{L^{\gamma(\rho)}(0,T; L^\rho)}
    \ls &C\norm{u_0}_\Sigma+CT^{1-\eps-\frac{2}{\gamma(\rho)}}\norm{u}_{L^\infty(0,T; H^1)}^{2\sigma}
    \nnorm{u}_{L^{\gamma(\rho)}(0,T; L^\rho)}\no\\
    \ls &M/2+CT^{1-\eps-\frac{2}{\gamma(\rho)}}M^{2\sigma}M\no\\
    \ls &M,
\end{align}
where we have taken $T\in (0,\pi/2\omega]$ so small that
$CT^{1-\eps-\frac{2}{\gamma(\rho)}}M^{2\sigma}\ls 1/2$. Similar to
the above, a straightforward computation shows that it holds
\begin{align}\label{non2}
    d(\T u,\T v)
    \ls &CT^{1-\eps-\frac{2}{\gamma(\rho)}}\left(\norm{u}_{L^\infty(0,T; H^1)}^{2\sigma}
    +\norm{v}_{L^\infty(0,T; H^1)}^{2\sigma}\right)
    \nnorm{u-v}_{L^{\gamma(\rho)}(0,T; L^\rho)}\no\\
    \ls &CT^{1-\eps-\frac{2}{\gamma(\rho)}}M^{2\sigma}d(u,v)\no\\
    \ls &\frac{1}{2}d(u,v).
\end{align}
Hence, $\T$ is a contracted mapping from the Banach space
$(\mathcal{D},d)$ to itself. By the Banach contraction mapping
principle, we know that there exists a unique solution
$u\in(\mathcal{D},d)$ to \eqref{gpe}--\eqref{data}. In view of the
conservation laws, we can use the standard argument to extend it
uniquely to a solution at the interval $[0,\pi/2\omega]$ which
satisfies for any $t\in[0,\pi/2\omega]$ and $\rho\in[2,\infty)$
\begin{align*}
    u(t,x), J(t)u(t,x), H(t)u(t,x)\in \mathcal{C}(0,\pi/2\omega;
    L^2(\Real^2))\cap L^{\gamma(\rho)}(0,\pi/2\omega; L^\rho(\Real^2)).
\end{align*}

Then, we can extend the above solution to a global one by
translation. In fact, in order to get the solution in the interval
$(\pi/2\omega, \pi/\omega]$, we can apply a translation
transformation with respect to the time variable $t$ such that the
initial data $u(\pi/2\omega)$ are replaced by $\tilde{u}(0)$. Let
$\tilde{u}(t,x):=u(t-\pi/2\omega,x)$, then we have from the original
equation with initial data $u(\pi/2\omega,x)$
\begin{align}
    &i\tilde{u}_t+\frac{1}{2}\Delta \tilde{u}=\frac{\omega^2}{2}\abs{x}^2 \tilde{u}+\beta\abs{\tilde{u}}^{2\sigma}\tilde{u}
    -\omega L_z \tilde{u},\;
    x\in\Real^2,\, t\gs 0, \label{gpe2}\\
    &\tilde{u}(0,x)=\tilde{u}_0(x):=u(\pi/2\omega,x),\; x\in\Real^2.\label{data2}
\end{align}
In the same way, we can get a solution $\tilde{u}(t,x)$  of
\eqref{gpe2}--\eqref{data2} for $t\in[0,\pi/2\omega]$. It is also a
solution $u(t,x)$ of \eqref{gpe}--\eqref{data} for
$t\in[\pi/2\omega, \pi/\omega]$ and it is unique. Thus, by an
induction argument with the help of those conserved identities
stated in Proposition~\ref{pr.cons}, we can obtain a global solution
$u(t,x)$ of \eqref{gpe}--\eqref{data} satisfying for any
$T\in(0,\infty)$
\begin{align*}
    u(t,x), J(t)u(t,x), H(t)u(t,x)\in \mathcal{C}(\Real;
    L^2(\Real^2))\cap L^{\gamma(\rho)}(0,T; L^\rho(\Real^2)).
\end{align*}

Therefore, we have completed the proof of the main theorem.

\section*{Acknowledgments}

C.C.Hao was partially supported by the Scientific Research Startup
Special Foundation on Excellent PhD Thesis and Presidential Award of
Chinese Academy of Sciences, NSFC (Grant No. 10601061), and the
Innovation Funds of AMSS, CAS of China. L.Hsiao was partially
supported by NSFC (Grant No. 10431060). H.L.Li was partially
supported by NSFC (Grant No. 10431060) and Beijing Nova Program.\\


\end{document}